\documentclass[a4paper,11pt]{article}
\usepackage{paralist}   

\usepackage{amsmath}
\usepackage{amssymb}
\usepackage{amsfonts}
\usepackage{color}

\usepackage{verbatim}

\usepackage{graphicx}
\usepackage{subfig}
\usepackage{bm}
\usepackage{tikz}
\newcommand{\red}[1]{{\color{red} #1 }}

\usepackage[bookmarks,bookmarksopen=false,
                  pdfauthor={Autor},%
                  pdftitle={Titel},%
                  colorlinks=true,%
                  linkcolor=blue,%
                  citecolor=blue,%
                  urlcolor=blue,%
                ]{hyperref}

\newtheorem{theorem}{Theorem}
\newtheorem{definition}{Definition}
\newtheorem{proposition}[theorem]{Proposition}
\newtheorem{corollary}[theorem]{Corollary}
\newtheorem{conjecture}{Conjecture}
\newtheorem{lemma}[theorem]{Lemma}
\newtheorem{remark}[theorem]{Remark}
\newtheorem{example}{Example}
\DeclareMathOperator{\esssup}{ess\,sup}
\DeclareMathOperator{\essinf}{ess\,inf}
\newcommand{\proof}{\textbf{Proof:\ }}
\newcommand{\pbox}{\hfill$\Box$\\}

\newcommand{\R}{\mathbb{R}}
\newcommand{\N}{\mathbb{N}}
\newcommand{\C}{\mathbb{C}}
\newcommand{\Z}{\mathbb{Z}}

\newcommand{\I}{\mathcal{I}}
\renewcommand{\H}{\mathcal{H}}
\newcommand{\dom}{{\sf Dom}\,}
\newcommand{\ran}{{\sf Ran}\,}
\newcommand{\Ker}{{\sf Ker}\,}
\renewcommand{\ker}{{\sf Ker}\,}

\newcommand{\norm}[2]{\left\| #2 \right\|_{#1}}
\newcommand{\ip}[2]{ \langle {#1} ,{#2}  \rangle}
\def\V{{\mathcal V}}
\definecolor{darkviolet}{rgb}{0.58,0,0.83} 
\newcommand{\xxl}[1]{{\color{darkviolet} #1 }}

\begin{document}
\title{Frames, their relatives and reproducing kernel Hilbert spaces}

\date{}
\author{\vspace{0.3cm} Michael Speckbacher\footnote{speckbacher@kfs.oeaw.ac.at}\hspace{0.15cm} and Peter Balazs\footnote{peter.balazs@oeaw.ac.at}\\ 
Acoustics Research Institute, Austrian Academy of Sciences,\\
Wohllebengasse 12-14, 1040 Vienna, Austria}
\maketitle

\begin{abstract}\noindent
This paper considers different facets of the interplay between reproducing kernel Hilbert spaces (RKHS) and stable analysis/synthesis processes: First, we analyze the structure of the reproducing kernel of a RKHS using frames and
reproducing pairs.
Second, we  present a new approach to prove the result that finite redundancy of a continuous frame implies atomic structure of the underlying measure space. Our proof uses the RKHS structure of the range of the analysis operator. 
This in turn implies that all the attempts to extend the notion of Riesz basis to general measure spaces are fruitless since every such family can be identified with a discrete Riesz basis.  
Finally, we show how the range of the analysis operators of a reproducing pair can be equipped with a RKHS structure.
\end{abstract}


\section{Introduction}
Reproducing kernel Hilbert spaces (RKHS) were introduced by Zaremba \cite{za1907} and Mercer \cite{me1909} and systematically studied by Aronszajn \cite{aro50}  in 1950. These spaces  play an important role in many diverse branches of mathematics such as complex analysis \cite{duschu04}, 
and  learning theory  \cite{schhesm01}. Another field with manifold  connections to RKHS, which is the main focus of this paper, is frame theory and related concepts.

 A mapping $x\mapsto\Psi_x\in \H$ is called a \emph{continuous frame} if there exist constants $m,M>0$ such that 
\begin{equation}\label{intro-frame}
m\|f\|^2\leq\int_X|\langle f,\Psi_x\rangle|^2d\mu(x)\leq M\|f\|^2,\qquad \forall
f\in\H.
\end{equation}
Frames have proven to be a viable tool in many different fields such as signal processing \cite{nsdgt10} 
acoustics \cite{framepsycho16} 
 or mathematical physics \cite{alanga00,xxlbayasg11}. It is however not always possible to satisfy both inequalities, which is why new concepts like semi-frames \cite{jpaxxl09,jpaxxl12} and reproducing pairs \cite{ansptr15,spexxl14,spexxl16} were introduced recently. An upper (resp. lower) semi-frame is a complete system that only satisfies the upper (resp. lower) frame inequality in \eqref{intro-frame}. A reproducing pair is a pair of mappings $(\Psi,\Phi)$  that generates a bounded and boundedly invertible analysis/synthesis process, i.e., 
$$
S_{\Psi,\Phi}f=\int_X\langle f,\Psi_x\rangle \Phi_x d\mu(x) \in  GL(\H),
$$
 without assuming any frame inequality for neither $\Psi$ nor $\Phi$.

This paper is divided into  three  main parts, portraying different facets of the interplay  between frames, reproducing pairs, and RKHS. 

In Section~\ref{sec:char-RKHS}, we 
study systems taking values in RKHS. In particular, we give an explicit expression for the reproducing kernel in terms of a reproducing pair in Theorem~\ref{kernel-and-rep-pair} that extends the results from \cite{paul09,raca05} and introduce a  necessary condition for a family of vectors to form a frame, see Proposition~\ref{bessel-rkhs}.

In Section~\ref{sec:redundancy-section} we study the dependence of the redundancy of (semi-)frames on the structure of the measure space. In the discrete case, the redundancy of a frame measures, loosely speaking, how much a Hilbert space is oversampled by that frame, see for example \cite{bocaku11,cacahe11}. It is however impossible to directly translate this concept to continuous (semi-)frames. In \cite{hedera00} the authors therefore used a property of Riesz bases 
to define redundancy. A Riesz basis  is a discrete, non-redundant frame, i.e., its analysis operator $C_\Psi:\H\rightarrow\ell^2$ is surjective. 
We therefore define  redundancy of a (semi-)frame $\Psi$  by
\begin{equation*}
R(\Psi):=\dim(\ran C_\Psi {}^\bot).
\end{equation*}
It has been observed in several articles \cite{hedera00,hogira13,jale15} that $R(\Psi)$ depends on the underlying measure space $(X,\mu)$. In particular, if a (lower semi-)frame has finite redundancy, then it follows that $(X,\mu)$ is atomic. The proofs in the aforementioned papers all rely in one way or the other on the following argument: If the redundancy of a frame is zero (finite), then
$$
\inf\big\{\mu(A):\ A\mbox{ measurable and }\mu(A)>0\big\}=C>0,
$$
which implies that $(X,\mu)$ is atomic. We present a new proof in Section~\ref{subsec:frames-and-redund} using that $\ran C_\Psi$ is a RKHS, which, in our opinion, better explains the underlying structure of the problem. 
This result also shows that all efforts to generalize the notion of Riesz bases ($R(\Psi)=0$), see \cite{arkatota12,gaha03}, to the setting of  measure spaces that are not atomic, are essentially an attempt to square the circle. The measure space of a continuous Riesz basis is always atomic and every such system can therefore be written as a discrete family of vectors, see Corollary~\ref{cont-riesz-is-discrete}.

For upper semi-frames on the other hand, there is no connection between the structure of the measure space and the redundancy. In particular, we give a sufficient condition for the existence of upper semi-frames indexed by a non-atomic measure space with redundancy zero, see Proposition~\ref{upper-zero-redun}.

Finally, Section~\ref{sec:rep-pair-rkhs} is concerned with characterizing the ranges of the analysis operators of a reproducing pair. The omission of the frame inequalities causes the inconvenience that $\ran C_\Psi$ and $\ran C_\Phi$ need no longer be contained in $L^2(X,\mu)$. Therefore, in \cite{ansptr15}  a pair of Hilbert spaces, intrinsically generated by a reproducing pair, was introduced to study this problem. We  demonstrate that these spaces are actually RKHS and calculate the reproducing kernels.


 
\section{Preliminaries}\label{sec:prel-rkhs}

\subsection{Atomic and non-atomic measures}
Throughout this paper we assume that $(X,\mu)$ is a  measure space, where $\mu$ is a non-trivial, $\sigma$-finite, and positive measure.
A measurable set $A \subset X$ is called an \emph{atom} if $\mu(A)>0$, and for any measurable subset $B\subset A$, 
with $\mu(B)<\mu(A)$,
it holds $\mu(B)=0$.
A measure space is called  \emph{atomic} if there exists a partition $\{A_n\}_{n\in \N}$ of $X$
consisting of atoms and null sets.
The space $(X,\mu)$ is called  \emph{non-atomic} if there are no atoms in $(X,\mu)$. To our knowledge there is no term to denote a
measure space that is 
not atomic. In order to avoid any confusion with non-atomic spaces, we therefore call a measure space \emph{an-atomic} if it is not atomic.

A well-known result by Sierpi{\'n}ski states that non-atomic measures take a continuity of values.
\begin{theorem}[Sierpi{\'n}ski \cite{sie22}]\label{sierpinski-thm}
 Let $(X,\mu)$ be non-atomic and let $A\subset X$ be measurable with $\mu(A)>0$. 
 For every $0\leq b\leq \mu(A)$,
 there exists $B\subset A$ such that $\mu(B)=b$.
\end{theorem}
Since we could not find any reference for the second part of the following lemma, we  provide a proof in the appendix.
\begin{lemma}\label{not-atomic-non-atomic}
Let $(X,\mu)$ be a $\sigma$-finite measure space.\begin{enumerate}[(i)]\item There exists $\mu_a$ atomic and $\mu_c$ non-atomic such that 
\begin{equation}\label{measure-partition}
\mu=\mu_a+\mu_c.
\end{equation}
\item If $(X,\mu)$ is an-atomic, then there exists $A\subset X$ with $\mu(A)>0$ and $(A,\mu)$ non-atomic.
\end{enumerate}
\end{lemma}

\subsection{Continuous frames, semi-frames and reproducing pairs}

Let $\mathcal{H}$ be a separable Hilbert space. We denote by $GL(\H)$ the space of bounded linear operators on $\H$ with bounded inverse.
\begin{definition}\label{def-cont-frame}
A mapping $\Psi:X\rightarrow \mathcal{H}$  is called a continuous frame if
\begin{enumerate}[(i)]
 \item $\Psi$ is weakly measurable, that is, $x\mapsto\langle f,\Psi_x\rangle$ is a measurable function for every
 $f\in\mathcal{H}$, 
 \item there exist positive constants $m,M>0$ such that
 \begin{equation}\label{frame-condition}
  m\left\|f\right\|^2\leq\int_{X}\left|\langle f,\Psi_x\rangle\right|^2d\mu(x)\leq M\left\|f\right\|^2,\qquad
  \forall f\in\mathcal{H}.
 \end{equation}
\end{enumerate}
The mapping $\Psi$ is called complete, if for every $f\in\H\backslash\{0\}$
$$
0<\int_{X}\left|\langle f,\Psi_x\rangle\right|^2d\mu(x).
$$
\end{definition}
The constants $m,M$ are called the frame bounds and $\Psi$ is called  Bessel if  the second inequality in 
\eqref{frame-condition} is satisfied. If $m=M=1$, then $\Psi$ is called a Parseval frame. 
If  $(X,\mu)$ is a countable set equipped with
the counting measure, then one recovers the classical definition of a discrete frame, see for example \cite{christ1}. For a  self-contained introduction to continuous 
frames, we refer the reader to \cite{ranade06}.

The fundamental operators in frame theory are given by the  \emph{analysis operator}
$
 C_\Psi:\mathcal{H}\rightarrow L^2(X,\mu)$, $C_\Psi f(x):=\langle f,\Psi_x\rangle,
$
and the  \emph{synthesis operator}
\begin{equation*}
 D_\Psi:{ {{L}^2(X,\mu)}}\rightarrow \mathcal{H},\ \ \ D_\Psi F:=\int_X F (x)\Psi_x d\mu(x),
\end{equation*} 
where the integral is defined weakly. 
The 
\emph{frame operator} $S_\Psi\in GL(\H)$ is defined as the composition of $C_\Psi$ and $D_\Psi$
\begin{equation*}
 S_\Psi:\mathcal{H}\rightarrow \mathcal{H},\ \ \ 
 S_\Psi f:=D_\Psi C_\Psi f=\int_X\langle f,\Psi_x\rangle \Psi_xd\mu(x).
\end{equation*}
Every frame $\Phi$ satisfying
\begin{equation*}
f:=D_\Psi C_{\Phi} f=D_{\Phi}C_\Psi f,\qquad \forall f\in \mathcal{H},
\end{equation*}
is called a \emph{dual frame} for $\Psi$. 

There exists a variety of interesting complete systems that do not meet both frame conditions. Several concepts to generalize the frame property were thus introduced and studied.
An \emph{upper  semi-frame} is a complete Bessel system, i.e.,
$$
0<\int_{X}\left|\langle f,\Psi_x\rangle\right|^2d\mu(x)\leq M\|f\|^2,\qquad \forall f\in\H\backslash\{0\},
$$ whereas a \emph{lower semi-frame} satisfies
$$
m \|f\|^2\leq \int_{X}\left|\langle f,\Psi_x\rangle\right|^2d\mu(x),
$$ 
see  \cite{jpaxxl09,
jpaxxl12}.

Another generalization  is  the concept of reproducing pairs, defined in \cite{spexxl14} 
and further investigated in \cite{ansptr15,antra16,spexxl16}. Here, one considers a pair of mappings instead of a single one without any assumption on frame inequalities.
\begin{definition}\label{rep-pair-definition}
 Let  $\Psi,\Phi:X\rightarrow\mathcal{H}$ be weakly measurable.
 The pair $(\Psi,\Phi)$ is called a  reproducing pair for $\mathcal{H}$ if the  operator 
 $S_{\Psi,\Phi}:\mathcal{H}\rightarrow \mathcal{H}$, defined by
 \begin{equation}\label{rep-pair-def}
 \langle S_{\Psi,\Phi} f,g\rangle:=\int_X \langle f,\Psi_x\rangle \langle\Phi_x,g\rangle d\mu(x),
 \end{equation}
is an element of $GL(\mathcal{H})$.
\end{definition}


\subsection{Reproducing kernel Hilbert spaces (RKHS)}\label{subsec:prel-rkhs}
Let $\mathcal{F}(\Omega,\mathbb{C})$ denote the vector space of all functions $f:\Omega\rightarrow\mathbb{C}$. 

\begin{definition}
A Hilbert space $\H_{K}\subset \mathcal{F}(\Omega,\C)$  is called  a  reproducing kernel Hilbert space (RKHS), if 
the point evaluation functional $\delta_z:\H_{K}\rightarrow\C$,  $\delta_z(f):=f(z)$ is bounded for 
  every $z\in \Omega$, that is, if there exists $C_z>0$ such that $|\delta_z(f)|\leq C_z\|f\|$, for every $f\in\H_K$.
\end{definition}
As $\delta_z$ is  bounded, there exists a unique vector $k_z\in \H_{K}$ such that $f(z)=\langle f,k_z\rangle,$ for every $f\in\H_{K}$.
The function $K(z,w):=k_w(z)=\langle k_w,k_z\rangle$ is called the  reproducing kernel for $\H_{K}$. The 
reproducing kernel is unique, $K(z,w)=\overline{K(w,z)}$ and 
its diagonal reads
$$K(z,z)=\langle k_z,k_z\rangle=\|k_z\|^2=\sup\big\{|f(z)|^2:\ f\in\H_K,\ \|f\|=1\big\}.$$
The following result can be found in \cite[Theorem 3.1 and 3.2]{alanga93}.
\begin{theorem}\label{charact-of-RKHS}
 If $\H_{K}$ is a RKHS and $\{\phi_i\}_{i\in\I}\subset \H_{K}$ an orthonormal basis, then 
  \begin{equation}K(z,w)=\sum_{i\in\I}\phi_i(z)\overline{\phi_i(w)},  \end{equation}
  with pointwise convergence of the series. In particular, 
  \begin{equation}\label{pointwise-l2-onb}
0<\sum_{i\in\I}|\phi_i(z)|^2=K(z,z)<\infty,\qquad \forall z\in \Omega.
\end{equation}
\\
Conversely, if  there exists an orthonormal basis for a Hilbert space $\H_K\subset \mathcal{F}(\Omega,\C)$ that satisfies \eqref{pointwise-l2-onb},
then $\H_{K}$ can be identified with a RKHS consisting of functions $f:\Omega\rightarrow \C$.

\end{theorem}
For a thorough introduction to RKHS we refer the reader to \cite{aro50,paul09}. 

\section{Frames and reproducing pairs taking values in a RKHS}\label{sec:char-RKHS}
In this section, we  investigate the pointwise behavior of frames in RKHS, characterize the reproducing kernel and introduce sufficient conditions on a frame that ensures the existence of a reproducing kernel.\\
The following result adapts the  
arguments of \cite[Theorem 3.12]{paul09} to reproducing pairs.
\begin{theorem}\label{kernel-and-rep-pair}
 Let $\H_{K}$ be a RKHS and $\Psi=\{\phi_i\}_{i\in\I},\ \Phi=\{\psi_i\}_{i\in\I}\subset\H_K$. The pair $(\Psi,\Phi)$
 is a reproducing pair for $\H_K$ if and only if there exists $A\in GL(\H_K)$ such that
\begin{equation}\label{rep-prod-char-ker}
K(z,w)=\sum_{i\in\I} (A\phi_i)(z)\overline{\psi_i(w)}=\sum_{i\in\I} (A^\ast\psi_i)(z)\overline{\phi_i(w)},\qquad \forall z,w\in\Omega,
\end{equation}
 where the series converges pointwise. In particular, $A$ is uniquely given by $S_{\Psi,\Phi}^{-1}$.
 \end{theorem}
\proof If $(\Psi,\Phi)$ is a reproducing pair, then 
$$
K(z,w)=\ip{k_w}{k_z}=\sum_{i\in\I} \ip{k_w}{\psi_i}\ip{S_{\Psi,\Phi}^{-1}\phi_i}{k_z}=\sum_{i\in\I} \overline{\psi_i(w)}(S_{\Psi,\Phi}^{-1}\phi_i)(z).
$$
Conversely, assume that $K$ is given by \eqref{rep-prod-char-ker}. If $f,g\in \mbox{span}\{k_z: z\in \Omega\}$, that is, there exist 
$\alpha_n,\beta_m\in \C$, and $z_n,w_m\in\Omega$,
such that $f=\sum_n^N \alpha_n k_{z_n}$ and $g=\sum_m^N \beta_m k_{w_m}$, then
\begin{align*}
\ip{f}{g}&=\sum_{n,m=1}^N\alpha_n\overline{\beta_m}\ip{k_{z_n}}{k_{w_m}}=\sum_{n,m=1}^N\alpha_n\overline{\beta_m}K(w_m,z_n)
\\
&=\sum_{n,m=1}^N\alpha_n\overline{\beta_m}\sum_{i\in\I}(A\phi_i)(w_m)\overline{\psi_i(z_n)}
\\&=
\sum_{n,m=1}^N\alpha_n\overline{\beta_m}\sum_{i\in\I}\ip{k_{z_n}}{\psi_i}\ip{A\phi_i}{k_{w_m}}
\\
&=\sum_{i\in\I}\Big\langle\sum_{n=1}^N\alpha_n k_{z_n},\psi_i\Big\rangle\Big\langle A\phi_i,\sum_{m=1}^N\beta_mk_{w_m}\Big\rangle
\\
&=\sum_{i\in\I}\ip{f}{\psi_i}\ip{A\phi_i}{g}=\ip{AS_{\Psi,\Phi}f}{g}.
\end{align*}
In \cite[Proposition 3.1]{paul09} it is shown that $\mbox{span}\{k_z: z\in \Omega\}$ is dense in $\H_K$. Therefore, it follows
that $AS_{\Psi,\Phi}=I$. As $A\in GL(\H_K)$ we may conclude that $S_{\Psi,\Phi}\in GL(\H_K)$, that is,
$(\Psi,\Phi)$ is a reproducing pair.
\pbox
\begin{remark}
For certain special cases, Theorem~\ref{kernel-and-rep-pair} is already  known. In particular, the result can be found in  \cite[Theorem 7]{raca05} if $\Psi$ and $\Phi$ are  dual frames, and in \cite[Theorem 3.12]{paul09} if $\Psi=\Phi$ is a Parseval frame. 
\end{remark}


\begin{proposition}\label{bessel-rkhs}
Let $\H_K$ be a RKHS and $\{\psi_i\}_{i\in\I}\subset \H_K$. 
\begin{enumerate}[(i)]
\item \label{item1}
If the family $\{\psi_i\}_{i\in\I}$ is Bessel, then
\begin{equation}\label{eq-lower-rkhs}
\frac{\sum_{i\in\I}|\psi_i(z)|^2}{K(z,z)}\leq M,\qquad \forall\ z\in \Omega.
\end{equation}
\item \label{item2}
If $\{\psi_i\}_{i\in\I}$ satisfies the lower frame inequality, then 
\begin{equation}\label{eq-lower-rkhs2}
0<m\leq\frac{\sum_{i\in\I}|\psi_i(z)|^2}{K(z,z)},\qquad \forall\ z\in \Omega.
\end{equation}
\end{enumerate}
\end{proposition}
\proof If $\{\psi_i\}_{i\in\I}$ is Bessel, then, for every $z\in \Omega$, it holds
$$ \sum_{i\in\I}|\psi_i(z)|^2=\sum_{i\in\I}|\ip{k_z}{\psi_i}|^2\leq M \|k_z\|^2=M K(z,z).$$
The same argument shows the lower bound in \eqref{eq-lower-rkhs2} if $\{\psi_i\}_{i\in\I}$ is a lower semi-frame.  \pbox

\begin{remark}\label{discrete-subspace-rkhs}
 (i) 
{The converse statements of Proposition~\ref{bessel-rkhs} are in general false.
First, consider the sequence $\psi_1=(1,0,\ldots),\ \psi_2=(0,1,1,0,\ldots),\psi_3=(0,0,0,1,1,1,0\ldots),\ldots\in\H_K=\ell(\N)$, i.e. $\psi_n(l)=1,$ if $l\in n(n-1)/2+1,\ldots,n(n+1)/2$, and zero otherwise. This system satisfies \eqref{eq-lower-rkhs} with $M=1$ but $\{\psi_n\}_{n\in\N}$ is not Bessel. To see this, recall that a sequence is Bessel if and only if its synthesis operator $D_\Psi$ is bounded. Taking $a=\{1/n\}_{n\in\N}$ gives $$\|D_\Psi a\|_{\ell^2(\N)}^2=\left\|\left(1,\frac{1}{2},\frac{1}{2},\frac{1}{3},\frac{1}{3},\frac{1}{3},\ldots\right)\right\|_{\ell^2(\N)}^2=\sum_{k\in\N}\frac{1}{k}=\infty.$$
}
  To see that the reverse of \eqref{eq-lower-rkhs2} is false, consider $\H_K=\C^2$, $\Omega=\{1,2\}$ and $\psi=(1,1)$. Then $K(z,z)=1$ and $|\psi(1)|^2=|\psi(2)|^2=1$, but $\psi$ is not spanning. 
 
The reverse statement remains false if one additionally assumes that $\{\psi_i\}_{i\in\I}$ is complete. Let us construct a counterexample from the system of integer time-frequency shifts of the Gaussian window $\varphi (t) = 2^{1/4} e^{- \pi t^2}$, which is a complete Bessel sequence (i.e. an upper semi-frame) but not a frame in $L^2(\R)$, see \cite{groe1,spexxl16}. Let $z=(x,\omega)\in \R^2$, and $\pi(z)f(t):=f(t-x)e^{2\pi i\omega t}$. The short-time Fourier transform with window $\varphi$ is defined as
 $
 V_\varphi f(z)=\langle f,\pi(z)\varphi\rangle ,
 $
 and the space 
 $$\H_K:=\{V_\varphi f:\ f\in L^2(\R)\}\subset L^2(\R^2),$$ forms a RKHS with kernel $K(z,w)=\langle \pi(w)\varphi,\pi(z)\varphi\rangle$. If  $F=V_\varphi f\in \H_K$ for some $f\in L^2(\R)$, then 
 $$
\sum_{\lambda\in\Z} |\langle F,K(\cdot,\lambda)\rangle|^2=\sum_{\lambda\in\Z} |\langle f,\pi(\lambda)\varphi\rangle|^2,
 $$ 
 and the isometry property $\|F\|=\|f\|$ 
implies that $\{\psi_\lambda\}_{\lambda\in\Z^2}=\{K(\hspace{0.05cm} \cdot\hspace{0.05cm},\lambda)\}_{\lambda\in\Z^2}$ is a complete Bessel sequence for $\H_K$, but not a frame. Moreover, we have $|K(z,w)|=e^{-\pi|z-w|^2/2}$ by \cite[Lemma 1.5.2]{groe1}, which gives that
 $$
 z\mapsto \sum_{\lambda\in\Z^2}|\psi_\lambda(z)|^2=\sum_{\lambda\in\Z^2}e^{-\pi|z-\lambda|^2}>e^{-\pi/2},
 $$
 is a strictly positive, continuous, and $\Z^2$-periodic function. Therefore, as $|K(z,z)|=1$, \eqref{eq-lower-rkhs2} is satisfied with $m\geq e^{-\pi/2}$.

\medskip
\noindent(ii) Another consequence of Proposition~\ref{bessel-rkhs} is that, if $\{\psi_i\}_{i\in\I}\subset\H_K$ is a Bessel sequence, then
\begin{equation}\label{synthesis-like}
\Big\|\sum_{i\in\I}\psi_i(z)\psi_i\Big\|\leq\sqrt{M\cdot K(z,z)}<\infty, \qquad \forall\  z\in \Omega.
\end{equation}
\end{remark}
\medskip
Let us now consider the converse part of Theorem~\ref{charact-of-RKHS}, which gives a sufficient condition on orthonormal bases that ensures that the space is a RKHS. We generalize the result to a condition on discrete frames.
\begin{theorem}\label{charact-of-RKHS_2}
Let $\H_K\subset \mathcal{F}(X,\C)$ be a Hilbert space. If there exists a discrete frame $\Psi$ for $\H_K$ such that 
  \begin{equation}\label{pointwise-l2-weak}
0<\sum_{i\in\I} |\psi_i(z)|^2=C_z< \infty,\qquad \forall\ z\in \Omega,
\end{equation}
then $\H_{K}$ is a RKHS.
\end{theorem}
\proof If $\Psi$ is a frame such that \eqref{pointwise-l2-weak} holds, then 
\begin{equation}\label{eq-some-sum}
\sum_{i\in\I} \overline{\psi_i(z)}S_\Psi^{-1}\psi_i
\end{equation}
is well-defined in $\H_K$ for every $z\in \Omega$. Moreover, for $f\in\H_K$, we have 
\begin{align*}
|f(z)|&=\Big|\sum_{i\in\I} \ip{f}{S_\Psi^{-1}\psi_i}\psi_i(z)\Big|=
\Big|\Big\langle f,\sum_{i\in\I}\overline{\psi_i(z)} S_\Psi^{-1}\psi_i\Big\rangle\Big|\\ &\leq
 \|f\|\Big\|\sum_{i\in\I}\overline{\psi_i(z)} S_\Psi^{-1}\psi_i\Big\|\leq
{\frac{1}{m}}\|f\|\Big\|\sum_{i\in\I}\overline{\psi_i(z)} \psi_i\Big\|\\ 
&\leq {\frac{\sqrt{M}}{m}}\left(\sum_{i\in\I}|\psi_i(z)|^2\right)^{1/2}\|f\|={\frac{\sqrt{C_z M}}{m}} \|f\|.
\end{align*}
Hence, point evaluation is continuous.\pbox



\noindent Note that  all results of this section can be reformulated for continuous frames and reproducing pairs, if one replaces the index set $\I$ by $(X,\mu)$. 

\section{Continuous (semi-)frames and their redundancy}\label{sec:redundancy-section}
Redundancy of a discrete frame measures, roughly speaking, how much the frame is oversampling the Hilbert space $\H$.
A non-redundant discrete frame is in fact a Riesz basis, that is, a frame with the additional
property that $\ran C_\Psi=(\ker D_\Psi)^\bot=\ell^2(\I)$. This justifies to define redundancy for general families $\Psi$ by
\begin{equation}
 R(\Psi):=\dim(\ran C_\Psi{}^\bot).
\end{equation}
Note that $C_\Psi$ is a Fredholm operator if $\Psi$ is Bessel and has finite redundancy. 
In that case, $R(\Psi)=-\mbox{ind}(C_\Psi)=\mbox{ind}(D_\Psi)$, where $\mbox{ind}(A)$ denotes the Fredholm index of a bounded operator $A$.

\begin{example}\label{ex:redundancy}
Let $\{e_n\}_{n\in\N}$ be 
an orthonormal basis. If we define $\Psi=\{e_1\}\cup\{e_n\}_{n\in\N}$, and 
$\Phi=\{e_n\}_{n\in\N}\cup\{e_n\}_{n\in\N}$, then $R(\Psi)=1$ and $R(\Phi)=\infty$. 

 If $\Psi=\{f_n\}_{n=1}^N$ is a finite frame for $\C^d$, then $R(\Psi)=N-d$, whereas the classical definition of redundancy for finite frames gives $N/d$. 
\end{example}

\subsection{Frames and lower semi-frames}\label{subsec:frames-and-redund}

The main goal of this section is to give a new proof for Theorem~\ref{reproduced-result}, which connects finite redundancy with the structure of the measure space, a result that has already been stated in \cite[Theorem 2]{hogira13}, \cite[Theorem 2.2]{hedera00} and \cite[Proposition 3.3]{jale15}.
 We thereby
use the property that $\ran C_\Psi$ forms a RKHS. In our opinion, this proof  does a better job explaining the inherent structure of continuous frames.

\begin{theorem}\label{reproduced-result}
 If a (lower semi-)frame $\Psi$ 
 has finite redundancy $ R(\Psi)<\infty$, then the measure space $(X,\mu)$ is atomic.
\end{theorem}
The converse is obviously not true, take for example the discrete family $\Phi$ in Example~\ref{ex:redundancy}.
\begin{corollary}
If $(X,\mu)$ is an-atomic and $\Psi$ a frame, then $R(\Psi)=\infty$ and there exist infinitely many dual frames for $\Psi$. 
\end{corollary}
{\proof Take an arbitrary pointwise defined function $\theta\in \ran C_\Psi^\bot$,  an arbitrary vector $g\in\H$ and define $\phi(x)=\theta(x)g$. Then $S_\Psi^{-1}\Psi+\Phi$ is a dual frame for $\Psi$.}
\pbox  

\noindent
We need to collect some auxiliary results in order to prove Theorem~\ref{reproduced-result}.

\begin{proposition}[\cite{gaha03}, Corollary 2.9]\label{exist-frame-for-RKHS}
 Let $\H_K$ be a  subspace of $L^2(X,\mu)$. If $\Psi$ satisfies the lower frame inequality, then $(\ran C_\Psi,\norm{2}{\cdot})$ is a RKHS. Moreover, the following are equivalent:
 \begin{enumerate}[(i)]
  \item There exists a continuous frame $\Psi$ such that $\ran(C_\Psi)=\H_K$,
  \item $\H_K$ is a RKHS.
 \end{enumerate}
\end{proposition}


\noindent In the following, we prove that $L^2(X,\mu)$ is not a RKHS if $(X,\mu)$ is an-atomic. This statement needs some elaboration to be precise. By definition, $L^2(X,\mu)$ is a space of equivalence classes of functions and does therefore not allow for pointwise evaluation. However, if we take an orthonormal basis $\{\phi_i\}_{i\in\I}$ for $L^2(X,\mu)$ and fix one particular representative for every basis element, then every vector $F\in L^2(X,\mu)$ can formally be written as
\begin{equation}\label{eq-f(x)}
F(x)=\sum_{i\in \I}\ip{F}{\phi_i}\phi_i(x),
\end{equation}
where the series converges for almost every $x\in X$. 

\begin{proposition} \label{L2-isnot-RKHS}
 If $(X,\mu)$ is a an-atomic measure space, then $L^2(X,\mu)$ is not a reproducing kernel Hilbert space.
\end{proposition}
\proof Every an-atomic measure space contains a non-atomic subspace by Lemma~\ref{not-atomic-non-atomic} $(ii)$. Without loss of generality we may therefore assume that $(X,\mu)$ is non-atomic.
Take $A\subset X$ with $\mu(A)>0$. We may assume without loss of generality that $\mu(A)=1$ as $(X,\mu)$ is $\sigma$-finite. Let $\{A_m\}_{m=1}^n$ be a partition of $A$ satisfying $\mu(A_m)=1/n,$ for every $m=1,...,n$. Such a partition exists by Theorem~\ref{sierpinski-thm}. 
Define $B^n\subset A$ by
$$
B^n:=\big\{x\in A:\ \chi_{A_m}(x)=1,\ \mbox{for some }m\in\{1,...,n\}\big\}. 
$$
Clearly, $\mu(B^n)=\mu(A)$ for every $n\in\N$. If we assume that $L^2(X,\mu)$ is a RKHS, then, for every $x\in B^n$, and some $m\in\{1,...,n\},$ one has $$|\chi_{A_{m}}(x)|^2=1=|\ip{\chi_{A_{m}}}{k_x}|^2\leq \|k_x\|^2/n.$$
In particular, $\|k_x\|^2\geq n$ for every $x\in B^n$. Setting $B:=\bigcap_{n\in\N}B^n$, one gets that $\mu(B)=\mu(A)$. Consequently, $\|k_x\|^2=K(x,x)=\infty$ for almost every $x\in A$, a contradiction to \eqref{pointwise-l2-onb}.\pbox

\begin{corollary} \label{noonbpointwise1}
 Let $(X,\mu)$ be an-atomic. There is no orthonormal basis $\{\phi_i\}_{i\in\I}\subset L^2(X,\mu)$, such that 
 $$
 \sum_{i\in\I}|\phi_i(x)|^2<\infty,\qquad \forall\ x\in X.
 $$
 In particular, for every orthonormal basis $\{\phi_i\}_{i\in\I}\subset L^2(X,\mu)$, there exists a set $A\subset X$
 of positive measure, such that 
  $$
 \sum_{i\in\I}|\phi_i(x)|^2=\infty,\qquad \forall x\in A.
 $$
\end{corollary}

\begin{proposition}\label{codim-of-RKHS-is-infty}
 If $(X,\mu)$ is an-atomic and $\H_K\subset L^2(X,\mu)$ is a RKHS, then $\dim (\H_K{}^\bot)=\infty$.
\end{proposition}
\proof If we assume that $\dim (\H_K{}^\bot)=N<\infty$, then any orthonormal basis $\{\phi_i\}_{i\in\I}$ of $\H_K$ can be complemented to a orthonormal basis of $L^2(X,\mu)$ using $N$ vectors $\{u_n\}_{n=1}^N$. In particular, we have by Theorem~\ref{charact-of-RKHS} that
$$
\sum_{i\in\I}|\phi_i(x)|^2+\sum_{n=1}^N|u_n(x)|^2<\infty,\qquad \forall\ x\in X,
$$
a contradiction to Corollary~\ref{noonbpointwise1}.\pbox


\noindent\textbf{Proof of Theorem~\ref{reproduced-result}:} The range of $C_\Psi$ is a RKHS by Propositions~\ref{exist-frame-for-RKHS}. 
 By  Proposition \ref{codim-of-RKHS-is-infty} it thus follows that either $R(\Psi)=\infty$ or $(X,\mu)$ is atomic. \hfill $\Box$

\begin{remark}
 The results of this section lead to interesting consequences in the context of quantum mechanics. Let us assume that $\Psi:X\rightarrow \H$ is a system of coherent states, see \cite{alanga00}. The probabilistic interpretation of quantum mechanics then states that the probability distribution of finding a system $f$ in the state $\Psi_x$ is given by $|\ip{f}{\Psi_x}|^2$. Hence, in light of Theorem~\ref{reproduced-result}, it follows that there is always a infinite dimensional subspace of  probability distributions that does not correspond to any physically feasible system.
\end{remark}



\subsection{Strictly continuous mappings}\label{subsec:strictly}
In the following, we  show that the discrete components of a continuous frame can be separated and we call the remaining remaining mapping a strictly continuous mapping.
\begin{definition}A mapping $\Psi:X\rightarrow \H$ is called strictly  continuous if $(X,\mu)$ is non-atomic and there exists no set $A\subset X$, $\mu(A)>0$, such that $C_\Psi f|_A$ is constant  for every $f\in\H$.
\end{definition}

\begin{example}
Let $G$ be a locally compact group, $\pi:G\rightarrow\mathcal{U}(\H)$ a square-integrable group representations for $G$, and  $\psi\in\H$ an admissible vector, see e.g. \cite{grmopa86} for more information. If the left Haar measure of $G$ is non-atomic, then $\Psi=\{\pi(g)\psi\}_{g\in G}$ is a strictly continuous frame, as $\langle f,\pi(g_1)\psi\rangle=\langle f,\pi(g_2)\psi\rangle$ for every $\forall f\in\H$, implies $g_1=g_2$. Short-time Fourier systems or continuous wavelet systems, see \cite{groe1}, are just two instances from this large class of strictly continuous mappings.
\end{example}

\noindent Throughout the rest of this section we show that every continuous frame can be decomposed into a discrete and a strictly continuous Bessel system.
\begin{lemma}[\cite{si96}, Theorem 3.8.1]\label{constant-on-atoms}
If  $A\subset X$ is an atom, and $F:X\rightarrow \C$ is measurable, then $F$ is constant
almost everywhere on $A$.
\end{lemma}

\begin{lemma}\label{atomic-means-discrete}
If $\Psi$ is Bessel, $A\subset X$ such that $\mu(A)>0$, and $\ip{f}{\Psi(\cdot)}$ is constant on $A$ for every $f\in \H$, then there exists a unique $\psi\in \H$ such that $$
\|C_\Psi f\|_2^2=\|C_\Psi f|_{X\backslash A}\|_2^2+|\ip{f}{\psi}|^2,\qquad \forall\ f\in\H.
$$
In particular, $\psi$ is weakly given by
\begin{equation}\label{defin-of-psi-const}
\ip{f}{\psi}:=\mu(A)^{-1/2}\int_{A}\ip{f}{\Psi_x}d\mu(x),\qquad \forall\ f\in\H.
\end{equation}
\end{lemma}
\proof First, observe that $\psi$ defined by \eqref{defin-of-psi-const} is unique for every $n\in\N$ by Riesz representation theorem 
$$
|\ip{f}{\psi}|\leq
\mu(A)^{-1/2}\int_{A}|\ip{f}{\Psi_x}|d\mu(x)
\leq \left(\int_{A}|\ip{f}{\Psi_x}|^2d\mu(x)\right)^{\frac{1}{2}}\leq \sqrt{M}\|f\|,
$$
where  $M$ is the upper frame bound of $\Psi$. Moreover, 
\begin{align*}
\int_X|\ip{f}{\Psi_x}|^2d\mu(x)&=\int_{X\backslash A}|\ip{f}{\Psi_x}|^2d\mu(x)+\int_{A}|\ip{f}{\Psi_x}|^2d\mu(x)
\\
&=\int_{X\backslash A}|\ip{f}{\Psi_x}|^2d\mu(x)+|\ip{f}{\psi}|^2
\end{align*}
where we used \eqref{defin-of-psi-const} and the fact that  $\ip{f}{\Psi(\cdot)}$ is almost everywhere constant on $A$.
\pbox

\begin{theorem}\label{decomp-of-frame}
Every frame $\Psi$ can be written as $\Psi=\Psi_d\cup \Psi_c$, where $\Psi_d$ is  a discrete Bessel system and $\Psi_c:(X_c,\mu_c)\rightarrow\H$ is a strictly continuous Bessel mapping with $X_c\subset X$.  
In particular, if $(X,\mu)$ is atomic, then $\Psi$ can be written as a discrete frame.
\end{theorem}
\proof By Lemma~\ref{not-atomic-non-atomic} $(i)$, every measure $\mu$ can be written as $\mu=\mu_a+\mu_c$, where $\mu_a$ is atomic and $\mu_c$ is non-atomic. By Lemma~\ref{constant-on-atoms} and \ref{atomic-means-discrete} we deduce that  $\Psi$ defined on $(X,\mu_a)$ can be identified with a discrete Bessel system $\Psi_a$. Let $X_d=\bigcup X_i\subset X$ be the disjoint union of all sets $X_i\subset X$, such that  $\mu_c(X_i)>0$, and $C_\Psi f|_{X_i}$ is constant for all $f\in\H$. Setting $A=X_i$ in \eqref{defin-of-psi-const}, then defines a family of vectors $\{\psi_i\}_{i\in\mathcal{I}}$. By definition $\Psi_c:=\Psi|_{X\backslash X_d}$ is a strictly continuous Bessel mapping. It therefore remains to show that $\mathcal{I}$ is countable.
This, however, is a direct consequence from the fact that $\sigma$-finite measure spaces can only be partitioned into countably many sets of positive measure. Hence,  $\Psi_d:=\Psi_a\cup \{\psi_i\}_{i\in\mathcal{I}}$ is a discrete Bessel sequence the result follows.\pbox

\noindent In an attempt to generalize the concept of Riesz bases, continuous Riesz bases \cite{arkatota12}  and Riesz-type mappings \cite{gaha03} were introduced.
It turns out that 
the these notions are equivalent and characterized as frames with  redundancy zero \cite[Proposition 2.5 \& Theorem 2.6]{arkatota12}.

\begin{corollary}\label{cont-riesz-is-discrete}
Every continuous Riesz basis (Riesz-type mapping) can be written as a discrete Riesz basis.
\end{corollary}
\proof If $\Psi$ is a continuous Riesz basis, then $R(\Psi)=0$ by definition. By Theorem~\ref{reproduced-result}, $(X,\mu)$ is 
atomic. Consequently, $\Psi$ corresponds to a discrete Riesz basis by  Theorem~\ref{decomp-of-frame}.\pbox

\subsection{Upper semi-frames}\label{sec:upper-semi}

Let us now illustrate how upper semi-frames behave fundamentally different than (lower semi-)frames. 
In particular, the closure of the range of the analysis operator is not necessarily a reproducing kernel Hilbert space and there exist upper semi-frames on non-atomic measure spaces with  redundancy zero (compare to Proposition~\ref{exist-frame-for-RKHS} and Theorem~\ref{reproduced-result}). 

\begin{example}\label{ex:affine}
In \cite{jpaxxl09,ansptr15} the following upper 
semi-frame has been studied. 
Take $\H_n:=L^2(\R^+,r^{n-1}d r)$, where $n\in\N$, and $(X,\mu)=(\R,d x)$.
 We use the following convention to denote the Fourier transform
 $$
 \widehat f(\omega)=\int_\R f(x)e^{-2\pi i x\omega}dx.
 $$
For $\psi\in \H_n$, we define the affine coherent state $\Psi_x$ by
$$
\Psi_x(r):=e^{-2\pi ixr}\psi(r),\qquad r\in\R^+,\ x\in\R.
$$
The mapping $\Psi$ forms an upper semi-frame if  
$\esssup_{r \in {\R}^{+}}r^{n-1}|\psi (r)|^{2}<\infty,$   
and $|\psi(r)|\neq 0$ for a.e. $r\in\R^+$.
The frame operator is then given by a multiplication operator on $\H_n$
$$
(S_\Psi f)(r)= r^{n-1}|\psi (r)|^{2} f(r).
$$
It is thus easy to see that $\Psi$ cannot form a frame since for every $\psi\in \H_n$,
$\essinf_{r \in {\R}^{+}} r^{n-1}|\psi (r)|^{2}=0$.
In \cite[Section 5.2]{ansptr15} it is shown that $\Ker D_\Psi=\mathcal F_+$, where 
$$\mathcal{F}_+ :=\{f\in L^2(\R):\ \widehat f(\omega)=0\ \text{for a.e. }\omega\geq0\}.$$
Clearly, $\overline{\ran C_\Psi}=(\ker D_\Psi)^\bot=\mathcal{F}_+{}^\bot=\mathcal{F}_-$, where
$$\mathcal{F}_- :=\{f\in L^2(\R):\ \widehat f(\omega)=0\ \text{for a.e. }\omega\leq0\}.$$
Therefore, $\Psi$ has infinite redundancy and a short argument shows that $\mathcal F_-$ is not a RKHS: 

 The dilation operator $D_a $, defined by $D_a f(x):=a^{-1/2}f(x/a)$, $a\in \R^+$,  acts isometrically on $\mathcal{F}_-$. Take $f\in\mathcal{F}_-$ with $\|f\|=1$ and $f(0)\neq0$, then $|D_af(0)|=|a^{-1/2}f(0)|\rightarrow\infty$, as $a\rightarrow 0$. Consequently, point evaluation cannot be continuous  and $\overline{\ran C_\Psi}=\mathcal{F}_-$ is not a RKHS.

The mapping $\Psi$ possesses several other interesting properties, see  \cite{ansptr15}. For instance, it forms
a total Bessel system with no dual, i.e.,  there is no mapping $\Phi$ such that $(\Psi,\Phi)$ generates a reproducing pair.
\end{example}

\begin{proposition}\label{upper-zero-redun}
If $(X,\mu)$ is a measure space, such that there exists an orthonormal basis $\{\psi_n\}_{n\in\N}$ of $L^2(X,\mu)$ satisfying
\begin{equation}\label{assumpt-on-onb}
\sup_{n\in\N}\sup_{x\in X}|\psi_n(x)|=C< \infty,
\end{equation}
then there exists an upper semi-frame $\Psi$ for $\H$ such that $\overline{\ran C_\Psi}=L^2(X,\mu)$.
In particular, $R(\Psi)=0$.
\end{proposition}
\proof Take  an arbitrary orthonormal basis $\{e_n\}_{n\in\N}$ of $\H$, and define
$$\Psi_x:=\sum_{n\in\N}n^{-1} e_n \psi_n(x).$$
The series converges absolutely in every point, and $\Psi$ is an upper semi-frame with the desired properties.
To see this, we first observe that $\Psi:X\rightarrow\H$ is well-defined 
as, for $x\in X$ fixed,
$$
|\ip{f}{\Psi_x}|\leq \sum_{n\in\N}|\ip{f}{e_n}n^{-1}\psi_n(x)|\leq \norm{}{f}\Big(\sum_{n\in\N}n^{-2}| \psi_n(x)|^2\Big)^{1/2}
$$
$$
\leq C\norm{}{f}\Big(\sum_{n\in\N}n^{-2}\Big)^{1/2}= \frac{\pi}{\sqrt{6}} C\norm{}{f},
$$
where we used \eqref{assumpt-on-onb} and Cauchy-Schwarz inequality.
Moreover,
$$
\int_X|\ip{f}{\Psi_x}|^2d\mu(x)\leq\int_X\norm{}{f}^2\sum_{n\in\N}n^{-2}| \psi_n(x)|^2d\mu(x)
$$
$$
=\norm{}{f}^2\sum_{n\in\N}n^{-2}\int_X| \psi_n(x)|^2 d\mu(x)=\norm{}{f}^2\sum_{n\in\N}n^{-2}=
\frac{\pi^2}{6}\norm{}{f}^2.
$$
 Since $\{\psi_n\}_{n\in\N}$ is an orthonormal basis of $L^2(X,\mu)$, it follows that $\Psi$ is total in $\H$, as 
$$
\int_X|\ip{f}{\Psi_x}|^2d\mu(x)= \int_X \sum_{n,k\in\N}\ip{f}{e_n}\ip{e_k}{f}(nk)^{-1}\psi_n(x)\overline{\psi_k(x)}d\mu(x)
$$
$$
=\sum_{n,k\in\N}\ip{f}{e_n}\ip{e_k}{f}(nk)^{-1}\delta_{n,k}=\sum_{n\in\N}|\ip{f}{e_n}|^2n^{-2}>0,
$$
for every $f\neq 0$.
Finally, the range of the analysis operator of the system $\{n^{-1} e_n \}_{n\in\N}$ is dense in $l^2(\N)$, which implies that $\ran C_\Psi$ is dense 
in $L^2(X,\mu)$.\pbox

\begin{example} Take the non-atomic measure space $(X,\mu)=(\mathbb{T},dx)$, where $\mathbb{T}$ denotes the torus, and the orthonormal Fourier basis 
$\psi_n(x)=e^{2\pi i xn},$ $n\in\Z$, then 
$$\sup_{n\in\Z}\sup_{x\in \mathbb{T}}|\psi_n(x)|=1.$$
Hence, there exists an upper semi-frame $\Psi$ with the closure of $\ran C_\Psi$ being $L^2(\mathbb T,dx)$, i.e., $R(\Psi)=0$.
\end{example}

\subsection{Existence of duals for lower semi-frames}\label{subsec:corr-proof}
In this section, we fill a gap in the arguments of the proof of \cite[Proposition 2.6]{jpaxxl09}, which states that for every lower semi-frame $\Psi$ there exists a dual Bessel mapping $\Phi$ such that $S_{\Psi,\Phi}=I$ on $\dom C_\Psi$. While the result itself is correct,
the construction of the dual system $\Phi$ in \cite{jpaxxl09}  is in general not  well-defined. In particular, $\Phi$ is 
defined there by 
$$
\Phi_x:=\sum_{n\in\N}\phi_n(x)V\phi_n=V\Big(\sum_{n\in\N}\phi_n(x)\phi_n\Big),
$$
where $V:L^2(X,\mu)\rightarrow \H$ is a bounded operator depending on $\Psi$ only and $\{\phi_n\}_{n\in\N}$ is an orthonormal basis
for $L^2(X,\mu)$. However, if $(X,\mu)$ is an-atomic, then by Corollary~\ref{noonbpointwise1} there exists a set of positive measure $A$ such that
$\sum|\phi_n(x)|^2=\infty,$  for all $x\in A$.
Thus,  $\Phi$ may be not well-defined on a set of positive measure.

\begin{proposition}[\cite{jpaxxl09}, Proposition 2.6]\label{corrected-result}
 If $\Psi$ is a lower semi-frame in $\H$, then there exists an upper semi-frame $\Phi$ such that 
 $$
 f=\int_X\ip{f}{\Psi_x}\Phi_xd\mu(x),\qquad \forall \ f\in \dom C_\Psi.
 $$
 Moreover, if $\dom C_\Psi$ is dense in $\H$, then 
  $$
 f=\int_X\ip{f}{\Phi_x}\Psi_xd\mu(x),\qquad \forall \ f\in \H.
 $$
\end{proposition}
\proof  If $\Psi$ is a lower semi-frame, then $\ran C_\Psi$ is a RKHS in $L^2(X,\mu)$ by Proposition~\ref{exist-frame-for-RKHS}.
Moreover, let $P$ denote the orthogonal projection from $L^2(X,\mu)$ onto $\ran C_\Psi$, and $\{e_n\}_{n\in\N}$ be an orthonormal
basis for $\H$.
Define the linear operator $V:L^2(X,\mu)\rightarrow \H$ by $V:=C_\Psi^{-1}$ on 
$\ran C_\Psi$ and $V:=0$ on $(\ran C_\Psi)^\bot$. Then $V$ is bounded and for all $f\in \dom C_\Psi,\ g\in\H$,
it holds 
\begin{align*}
\ip{f}{g}&=\ip{VC_\Psi f}{g}=\ip{C_\Psi f}{V^\ast g}_2\\ &=\Big\langle C_\Psi f,V^\ast \Big(\sum_{n\in\N}\ip{ g}{e_n}e_n\Big)\Big\rangle_2
=\Big\langle C_\Psi f,\sum_{n\in\N}\ip{g}{ e_n}V^\ast e_n\Big\rangle_2
\\ &=\Big\langle C_\Psi f,\sum_{n\in\N}\ip{g}{ e_n}P V^\ast e_n\Big\rangle_2
=\ip{C_\Psi f}{C_\Phi g}_2,
\end{align*}
where $\Phi_x:=\sum_{n\in\N}\overline{(PV^\ast e_n)}(x)e_n$. It remains to show that $\Phi_x$ is well-defined for every $x\in X$. Since $\{e_n\}_{n\in\N}$ is an orthonormal basis, if follows that $\Phi_x$ is well-defined
if and only if
$$
\sum_{n\in\N}|(PV^\ast e_n)(x)|^2<\infty,\qquad \forall\ x\in X.
$$ 
By Proposition~\ref{bessel-rkhs}, it is sufficient to show that $\{P V^\ast e_n\}_{n\in\N}$ is a Bessel sequence on $\ran C_\Psi$.
If $F\in \ran C_\Psi$, then
$$
\sum_{n\in\N}|\ip{F}{P V^\ast e_n}_2|^2=\sum_{n\in\N}|\ip{VPF}{e_n}|^2=\|VF\|^2\leq C\|F\|_2^2,
$$
as  $PF=F$ and $V$ is bounded. Finally, it remains to show that $\Phi$ is Bessel. If $f\in \H$, then
$$
\int_X|\ip{f}{\Phi_x}|^2d\mu(x)=\int_X\Big|\sum_{n\in\N}\ip{f}{e_n}PV^\ast e_n(x)\Big|^2d\mu(x)
$$
$$
=\big\|D_{PV^\ast e_n}\big(\{\ip{f}{e_n}\}_{n\in\N}\big)\big\|_2^2
\leq C\sum_{n\in\N}|\ip{f}{e_n}|^2=C\|f\|^2,
$$
as $\{P V^\ast e_n\}_{n\in\N}$ is Bessel.
\pbox

\begin{remark}
There is no analogue result of Proposition~\ref{corrected-result} if $\Psi$ is an upper semi-frame. In \cite{ansptr15} it is shown that the affine coherent state system presented in Section~\ref{sec:upper-semi} is a complete Bessel mapping with no dual.
\end{remark}

\section{Reproducing pairs and RKHSs}\label{sec:rep-pair-rkhs}

The absence of frame bounds causes trouble analyzing  $\ran C_\Psi$ and $\ran C_\Phi$ of a reproducing pair $(\Psi,\Phi)$. Without an upper frame bound it is no longer guaranteed that $\ran C_\Psi$ is a subspace of $L^2(X,\mu)$. The lower frame inequality, on the other hand, ensures that $\ran C_\Psi$ is a RKHS. In  \cite{ansptr15},
a construction of two mutually dual Hilbert spaces that are intrinsically generated by the pair $(\Psi,\Phi)$ is given.
Let us first recall some of the results before  explaining how  RKHS enter the picture.

Let  $\V_\Phi(X, \mu)$ be the space of all measurable functions  $F : X \to \C$ such that
$$\label{eq-Vphi}
\left|  \int_X  F(x)  \ip{\Phi_x}{g} d\mu(x) \right| \leq M \norm{}{g}, \qquad \forall\, g \in \H.
$$
Note that in general neither $\V_\Phi(X, \mu)\subset L^2(X, \mu)$ nor $ L^2(X, \mu)\subset \V_\Phi(X, \mu)$.
The linear map $T_\Phi :\V_\Phi(X, \mu) \rightarrow \H$  given weakly  by
\begin{equation}\label{def-T-phi}
\ip{T_\Phi   F}{g}  =\int_X  F(x)  \ip{\Phi_x}{g} d\mu(x) ,\qquad g\in\H,
\end{equation}
is thus well-defined by Riesz representation theorem, and can be seen as the natural extension of the synthesis operator $D_\Phi$ (defined on $\dom D_\Phi\subseteq L^2(X,\mu)$) to $\V_\Phi(X, \mu)$. 

Let $(\Psi,\Phi)$ be a reproducing pair. According to \cite{ansptr15},
it then holds that
\begin{equation}\label{dir-sum-rp}
\V_\Phi(X,\mu)=\ran C_\Psi\oplus \ker T_\Phi.
\end{equation}
This observation, together with the fact that $T_\Phi$ is general not one-to-one, motivates to define the redundancy for arbitrary complete mappings via 
\begin{equation}\label{redund-rep-pair}
R(\Phi):=\dim(\ker T_\Phi).
\end{equation}
We expect that similar results as in Section~\ref{subsec:frames-and-redund} hold for this definition.
\begin{conjecture} If $R(\Phi)<\infty$, then $(X,\mu)$ is atomic.
\end{conjecture}
The main difficulty is that there is no characterization of $\V_\Phi(X,\mu)$ that would allow  to treat the problem in a similar manner than in Section~\ref{subsec:frames-and-redund} using \eqref{dir-sum-rp}. It is in particular not even clear if $\V_\Phi(X,\mu)$ is normable.

Let us introduce the following vector space
$$ 
V_\Phi(X, \mu)= \V_\Phi(X, \mu)/{\Ker}\,T_\Phi,
$$
equipped with the inner product
$$
\ip{F}{G}_{\Phi}: =\ip{T_\phi F}{T_\phi G}, \ \  F,
G \in V_\Phi(X,\mu).
$$
This is indeed an inner product as $\ip{F}{F}_{\Phi}=0$ if and only if $F\in\Ker T_\Phi$.  Hence, $V_\Phi(X,\mu)$ forms a pre-Hilbert space and $T_\Phi:V_\Phi(X,\mu)\rightarrow \H$ is an isometry.
By \eqref{def-T-phi} $ \ip{\cdot}{\cdot}_{\Phi}$ can be written explicitly as
\begin{equation}\label{phi-inner-expl}
 \ip{F}{G}_{\Phi}=\int_X\int_X  F(x)  \ip{\Phi_x}{\Phi_y}\overline{G(y)} d\mu(x)d\mu(y).
\end{equation}

\noindent The following result answers the question if, given a mapping $\Phi$, 
there exist another mapping $\Psi$ such that $(\Psi,\Phi)$ forms a reproducing pair. 
 \begin{theorem}[\cite{ansptr15}, Theorem 4.1]\label{theo-partner}
Let $\Phi:X\rightarrow\H$ be a weakly measurable mapping and $\{e_i\}_{i\in\I}$ an orthonormal 
basis of $\H$. There exists another family $\Psi$, such that $(\Psi,\Phi)$ is a reproducing pair if and only if
\begin{enumerate}[(i)]\item $\ran  T_\phi =\H$, 
\item there exists $\{\mathcal{E}_i\}_{i\in\I}\subset \V_\Phi(X,\mu)$ satisfying $T_\Phi \mathcal{E}_i= e_i,\ \forall\ i\in\I,$ and
\begin{equation}\label{second-assumption}
\sum_{i\in\I}|\mathcal{E}_i(x)|^2<\infty,\qquad \forall\ x\in X.
\end{equation}
\end{enumerate}
 A reproducing partner $\Psi$ is then given by
\begin{equation}\label{def-repr-partn}
 \Psi_x:=\sum_{i\in\I}\overline{\mathcal{E}_i(x)}e_i.
\end{equation}
  \end{theorem}
 Theorem~\ref{theo-partner} is a powerful tool for the study of complete systems. It has for example been used to construct a reproducing partner for the Gabor system of integer time-frequency shifts of the Gaussian window  \cite{spexxl16} and to prove the non-existence of a dual for the system of affine coherent states in Example~\ref{ex:affine}.

Let us briefly discuss the conditions $(i)$ and $(ii)$ and some conceptual interpretations. For a complete system one can show that (under  mild conditions \cite[Lemma 2.2]{jpaxxl09}) $\overline{\ran D_\Phi}=\H$. It might therefore seem that $(i)$ is trivially satisfied for complete systems since $T_\Phi$ extends  $D_\Phi$ to its domain $\V_\Phi(X,\mu)$.  The complete upper semi-frame {from Example~\ref{ex:affine}} however does not satisfy $(i)$, see \cite[Section 6.2.3]{ansptr15}. 

\medskip

\noindent {\bf Coefficient map interpretation:} Property $(i)$ ensures the existence of a linear coefficient
map $A:\H\rightarrow \V_\Phi(X,\mu)$ satisfying $f=T_\Phi A(f)$ for every $f\in\H$. 
Property $(ii)$ then guarantees that  $A(f)$ can be calculated taking inner products of $f$ with a second mapping
$\Psi:X\rightarrow\H$. 

\medskip

\noindent {\bf RKHS interpretation:} If $(i)$ and $(ii)$ are satisfied, then it follows that
$\{\mathcal{E}_i\}_{i\in\I}$ forms an orthonormal family with respect to  $\ip{\cdot}{\cdot}_\Phi$, since  $(ii)$ implies
$$
\langle \mathcal{E}_i,\mathcal{E}_k\rangle_\Phi=\langle T_\Phi \mathcal{E}_i,T_\Phi \mathcal{E}_k\rangle=\langle e_i,e_k\rangle=\delta_{i,k}.
$$
Hence, $\{\mathcal{E}_i\}_{i\in\I}$ forms an orthonormal basis for
$$ 
\H_K^\Phi:=\overline{\mbox{span}\{\mathcal{E}_i:\ i\in\I\}}^{\|\cdot\|_\Phi}.
$$
Theorem~\ref{charact-of-RKHS} together with \eqref{second-assumption} thus ensure that $\H_K^\Phi$ is a RKHS.
Moreover, the definition of the reproducing partner $\Psi$ in \eqref{def-repr-partn} yields that
\begin{equation}\label{rkhs-ran-anal}\H_K^\Phi\simeq V_\Phi(X,\mu)\simeq(\ran C_\Psi,\|\cdot\|_\Phi).
\end{equation}
To put it another way, $(i)$ and $(ii)$ guarantee the existence of a  RKHS $\H_K^\Phi\subset \V_\Phi(X,\mu)$
reproducing $\H$ in the sense that $T_\Phi(\H_K^\Phi)=\H$.

\medskip
Let us assume that $(\Psi,\Phi)$ is a reproducing pair. There is a natural way to  generate frames on  $\H$ and $\H_K^\Phi$ via the analysis and synthesis operators.
\begin{proposition}\label{rep-pair-frame-dec}
Let $(\Psi,\Phi)$ be a reproducing pair for $\H$, $\{g_i\}_{i\in\I}$ a frame for $\H$, and $\{G_i\}_{i\in\I}$ a frame for 
$\H_K^\Phi$. 
If $H_i(x):=\ip{g_i}{\Psi_x}$ and $h_i:=T_\Phi G_i$, then
$\{H_i\}_{i\in\I}$ is a frame for $\H_K^\Phi$ and $\{h_i\}_{i\in\I}$ is a frame for $\H$.
\end{proposition}
\proof If $F\in \H_K^\Phi$, then 
$$
\sum_{i\in\I}|\langle F,H_i\rangle_\Phi|^2=\sum_{i\in\I}|\langle T_\Phi F,T_\Phi H_i\rangle|^2
=
\sum_{i\in\I}|\langle T_\Phi F,S_{\Psi,\Phi} g_i\rangle|^2
$$
$$
=\sum_{i\in\I}|\langle (S_{\Psi,\Phi})^\ast T_\Phi F, g_i\rangle|^2\leq M\|(S_{\Psi,\Phi})^\ast T_\Phi F\|^2
$$
$$
\leq M\|S_{\Psi,\Phi}\|^2\| T_\Phi F\|^2 =\widetilde M \|F\|_\Phi^2.
$$
The lower bound follows from the same argument as $(S_{\Psi,\Phi})^\ast$ is boundedly invertible.
Hence, $\{H_i\}_{i\in\I}$ is a frame for $\H_K^\Phi$. \\
For $f\in \H$, we have
$$
\|f\|=\|T_\Phi C_\Psi S_{\Psi,\Phi}^{-1}f\|=\|C_\Psi S_{\Psi,\Phi}^{-1}f\|_\Phi, 
$$
which, together with 
$$
\sum_{i\in\I}|\langle f,h_i\rangle|^2=\sum_{i\in\I}|\langle T_\Phi C_\Psi S_{\Psi,\Phi}^{-1}f,T_\Phi G_i\rangle|^2
=
\sum_{i\in\I}|\langle C_\Psi S_{\Psi,\Phi}^{-1}f,G_i\rangle_\Phi|^2,
$$
yields that $\{h_i\}_{i\in\I}$ is a frame for $\H$.
\pbox

\noindent The rest of this section is concerned with the explicit calculation of the reproducing kernel for $\H_K^\Phi$.
For   a reproducing pair $(\Psi,\Phi)$, there exists a similar characterization of the range of the 
 analysis operators as  for frames. 
In particular, if $R_{\Psi,\Phi}(x,y):=\ip{S_{\Psi,\Phi}^{-1}\Phi_y}{\Psi_x}$ defines the
 integral operator
$$
\mathcal{R}_{\Psi,\Phi}(F)(x):=\int_X F(y)R_{\Psi,\Phi}(x,y)d\mu(y),\qquad F\in \V_\Phi(X,\mu),
$$
then
by \cite[Proposition 2]{spexxl14} $\mathcal{R}_{\Psi,\Phi}(F)(x)=F(x)$ if and only if there 
exists $f\in \H$ such that $F(x)=\ip{f}{\Psi_x}$, for all $x\in X$.
However, $R_{\Psi,\Phi}$ is not the reproducing kernel for $\H_K^\Phi$ since the reproducing formula is based
on the inner product of $L^2(X,\mu)$ and not on $\langle \cdot,\cdot\rangle_\Phi$. 

If $F\in\ran C_\Psi$, then \eqref{phi-inner-expl} and the identity $f=T_\Phi C_\Psi S_{\Psi,\Phi}^{-1} f$ yield 
\begin{align*}
 F(x)&=\mathcal{R}_{\Psi,\Phi}(F)(x)=\int_X F(y)\ip{\Phi_y}{(S_{\Psi,\Phi}^{-1})^\ast\Psi_x}d\mu(y)\\
 &=\int_X \int_X F(y)\ip{\Phi_y}{\Phi_z}\ip{\Psi_z}{S_{\Psi,\Phi}^{-1}(S_{\Psi,\Phi}^{-1})^\ast\Psi_x}d\mu(z)d\mu(y)\\
\\
&=\Big\langle F,\big\langle (S_{\Psi,\Phi}^{-1})^\ast\Psi_x,(S_{\Psi,\Phi}^{-1})^\ast\Psi_\cdot\big\rangle\Big\rangle_\Phi
 \end{align*}
Hence, using $(S_{\Psi,\Phi}^{-1})^\ast=S_{\Phi,\Psi}^{-1}$, we finally obtain that
$$
K_{\Phi}(x,y)=\big\langle S_{\Phi,\Psi}^{-1}\Psi_x,S_{\Phi,\Psi}^{-1}\Psi_y\big\rangle
$$
is the reproducing kernel for $\H_K^\Phi$.

\section{Conclusion}
 With the results of Section~\ref{sec:redundancy-section} in mind, we suggest to use  the term 
\emph{continuous frame} only in the case of a strictly continuous frame, and 
\emph{semi-continuous} or \emph{semi-discrete} frame if it can be decomposed into nontrivial strictly continuous and  discrete parts. Moreover, the notion of a continuous Riesz basis (Riesz type mapping) should \textit{not} be used any further as every such system can be written as a discrete Riesz basis.


An interesting topic for future research is to find and study alternative notions of redundancy for continuous frames. A promising approach that may be adapted can be found in \cite{bacahela06}. Exploring the dependence of any notion of redundancy on the underlying measure space should remain a key task.


\section*{Appendix}

\textbf{Proof of Lemma~\ref{not-atomic-non-atomic}:} Ad $(i)$: See  \cite{fi72}.

Ad $(ii)$: Let $(X,\mu)$ be an-atomic. Let us assume to the contrary that 
for every measurable set $A\subset X$ with 
$\mu(A)>0$ there exists an atom $B\subset A$, and take $\{A_n\}_{n\in\mathbb{N}}\subset X$ to be a countable partition of $X$ by sets of finite measure. We  show that each $A_n$ can be partitioned into atoms and null sets, a contradiction. Assume without loss of generality that  $\mu(A_{1})>0$. By assumption, there exists an atom
$B_1\subset A_{1}$. 
If  $\mu(B_1)=\mu(A_{1})$, then $A_{1}$ is an atom. If  $0<\mu(B_1)<\mu(A_{1})$, then
$\mu(A_{1}\backslash B_1)>0$. Hence, there exists an atom $B_2\subset A_{1}\backslash B_1$ and the preceding
argument can be repeated. If  one has
$\mu\big(A_{1}\backslash \big(\bigcup_{k=1}^KB_k\big)\big)>0$ for all iteration steps $K$, then $\mu_K:=\mu\big(\bigcup_{k=1}^K B_k\big)$
defines a  strictly increasing sequence, bounded by $\mu(A_{1})$.
Hence, $\mu_K$ is convergent to some $\mu^\ast$ and the limit equals $\mu(A_{1})$. 
Indeed, if $\mu^\ast<\mu(A_1)$  then, by  assumption,
there exists an atom $B^\ast\subset A_{1}\backslash \bigcup_{k\in\N} B_k$ and
$$\mu\Big(\bigcup_{k\in\N} B_k\cup B^\ast\Big)>\mu^\ast,$$ 
a contradiction.
Consequently,
$
A_{1}=\bigcup_{k\in\N} B_k\cup N,$
where $N=A_{1}\backslash \bigcup_{k\in\N} B_k$ is of measure zero. In particular, we constructed a 
partition of $A_{1}$ consisting of atoms and null sets.
Repeating this argument for every $A_n$, $n\in\mathcal{I}$, with $\mu(A_n)>0$ shows that $(X,\mu)$ is atomic, a
contradiction. \pbox

\section*{Acknowledgement}
This work was funded by the Austrian Science Fund (FWF) START-project FLAME ('Frames and
Linear Operators for Acoustical Modeling and Parameter Estimation'; Y 551-N13).

The authors would like to thank Jean-Pierre Antoine for fruitful discussions on the physical interpretation
of the result on the redundancy of continuous frames. The second authors thanks Nora Simovich for help with typing.

\bibliographystyle{plain}
\bibliography{paperbib}
 \end{document}